\documentclass[12pt]{amsart}
\usepackage{amssymb}
\usepackage{calrsfs}
\usepackage{fullpage}

\usepackage{color}

\newcommand{\old}[1]{\relax} %old version
\newcommand{\ord}{\mathrm{ord}}

\newcommand{\be}{\begin{enumerate}}
\newcommand{\ee}{\end{enumerate}}

\newcommand{\calD}{{\mathcal D}}

\newcommand{\calP}{{\mathcal P}}

\newcommand{\calS}{{\mathcal S}}

\newcommand{\calW}{{\mathcal W}}

\newtheorem{introthm}{Theorem}

\newcommand{\C}{{\mathbb C}}

\newcommand{\Q}{{\mathbb Q}}

\newcommand{\Z}{{\mathbb Z}}

\newcommand{\pp}{{\mathfrak p}}

\newcommand{\qq}{{\mathfrak q}}

\newcommand{\Dd}{{\mathfrak D}}

\newtheorem{theorem}{Theorem}[section]
\newtheorem{lemma}[theorem]{Lemma}

\newtheorem{proposition}[theorem]{Proposition}

\theoremstyle{definition}
\newtheorem{definition}[theorem]{Definition}

\theoremstyle{remark}
\newtheorem{remark}[theorem]{Remark}

\newtheorem{notationassumption}[theorem]{Notation and Assumptions}

\newcounter{introrem}
\newtheorem{introremark}[introrem]{Remark}

\begin{document}
\bibliographystyle{plain}%
 \title{Diophantine Undecidability of Holomorphy Rings of Function Fields of Characteristic 0}
\date{\today}
\author{Laurent Moret-Bailly and Alexandra Shlapentokh}

\thanks{The second author has been partially supported by NSF grants DMS-0354907 and DMS-0650927.}
\address{IRMAR, Universit\'e de Rennes 1, Campus de Beaulieu, F-35042 Rennes
Cedex}\email{Laurent.Moret-Bailly@univ-rennes1.fr}\urladdr{perso.univ-rennes1.fr/laurent.moret-bailly/}\address{Department
of Mathematics \\ East Carolina University \\ Greenville, NC 27858}\email{shlapentokha@ecu.edu }
\urladdr{www.personal.ecu.edu/shlapentokha} 
\begin{abstract} Let
$K$ be a one-variable function field over a field of constants of characteristic 0. Let $R$ be a holomorphy subring of $K$, not equal to $K$. We prove the following undecidability results for $R$:
If $K$ is recursive, then Hilbert's Tenth Problem is undecidable in $R$. In general, there exist
$x_1,\ldots,x_n \in R$ such that there is no algorithm to tell whether a polynomial equation with
coefficients in $\Q(x_1,\ldots,x_n)$ has solutions in $R$.
 \end{abstract}%
\maketitle%

\section{Introduction}
The interest in the questions of existential definability and decidability over rings goes back to a question that
was posed by Hilbert: given an arbitrary polynomial equation in several variables over $\Z$, is there a uniform
algorithm to determine whether such an equation has solutions in $\Z$? This question, otherwise known as Hilbert's
Tenth Problem (``HTP'' in the future), has been answered negatively in the work of M. Davis, H. Putnam, J. Robinson
and Yu. Matijasevich. (See \cite{Da1},  \cite{Da2} or \cite{Mat}  for the details of the solution of the original
problem.)   Since the time when this result was obtained, similar questions have been raised for other fields and
rings. In other words, let $R$ be a  ring.  Then, given an arbitrary polynomial equation in several variables over
a recursive subring $R_0$ of  $R$, is there a uniform algorithm to determine whether such an equation has solutions
in $R$?  (If $R$ is countable and recursive then we can set $R_0=R$.)

 Depending on the nature of the ring the difficulty of answering the question can vary widely. By now,
a lot of work has been done to solve the problem over  some subrings of  number fields and  function fields,
including the fields themselves in the case of function fields.  However there remain quite a few open questions
which at the moment seem intractable.  Chief among these questions are arguably the  Diophantine
status of $\Q$ (and number fields in general), the rings of integers of an arbitrary number field, and
 an arbitrary function field of characteristic 0.

More details on the Diophantine problem over number fields and related issues can be found in \cite{CTSDS},
\cite{CPZ}, \cite{CZ}, \cite{Den1}, \cite{Den3}, \cite{Den2}, \cite{M1}, \cite{M2}, \cite{M3}, \cite{M4},
\cite{Ph1}, \cite{PO2}, \cite{PS}, \cite{Sha-Sh}, \cite{Sh2}, \cite{Sh5}, \cite{Sh17}, \cite{Sh1}, \cite{Sh3},
\cite{Sh21}, and  \cite{Sh26}. Results concerning function fields of positive characteristic can be found in
\cite{Den5}, \cite{Eis}, \cite{Ph3}, \cite{Ph7}, \cite{Sh30}, \cite{Sh13}, \cite{Sh24}, and  \cite{V}. Also, for a
general reference on the subject we suggest \cite{Lip4} and \cite{Sh34}.

It turned out that solving HTP over function fields of characteristic 0 was more difficult
  than over function fields of positive characteristic.  However we do know that HTP is
undecidable over many function fields  and rings of characteristic 0.  In particular,  we know that
HTP is undecidable over fields of functions of finite transcendence degree over constant fields that
are formally real or are subfields of finite extensions of $\Q_p$ for some odd rational prime $p$.
(These constant fields include number fields.)  Further, we also know that HTP is undecidable over function
fields of transcendence degree at least  2 whose field of constants is $\C$.  (See  \cite{Den4},
\cite{Eis2}, \cite{Eispadic}, \cite{K-R2}, \cite{K-R3}, \cite{MB3}, \cite{Z} for more details on
these field results).    We also have a few ring results: for rings of $\calS$-integers  and
semi-local rings over any field of constants, and some results for rings ``in the middle''.   (See
\cite{MB3}, \cite{Sh18}, \cite{Sh23}, \cite{Z2} for more details on ring results.)

One of the problems which was solved over global fields was the construction of an existential definition of order.
 In other words there exists an existential definition in the language of rings of the set of elements of a given
global field whose order at a fixed non-archimedean valuation is non-negative.  Over function fields of
characteristic 0 this was done successfully over a limited class of fields and the success depended heavily on the
nature of the field of constants.  As it turned out,  an existential definition of order was one of the two
ingredients  used for showing the Diophantine undecidability of one variable function fields of
characteristic 0.  The other ingredient was an elliptic curve of rank 1.  This plan was first
implemented by Denef for the formally real rational function fields in \cite{Den4}.   As will be
described below, the issue of finding the right elliptic curve has been solved in the greatest
possible generality by the first author in \cite{MB3}, but the issue of defining the order remains a
stumbling block.  So we can solve HTP precisely over those function fields in one variable where we
can define the order.  (See Corollary 10.3.3 of \cite{MB3}.)

To make the matters even more vexing, it is not hard to see that a definition of order together with
the Diophantine undecidability of any semilocal subring of a field implies Diophantine
undecidability of any ring ``in the middle'', i.e. any holomorphy ring, pretty much in the same
fashion as the Diophantine undecidability of a domain follows from the Diophantine undecidability of
its field of fractions. (Both cases also require being able to define the set of non-zero
elements. Fortunately, we know how to do it in all cases of interest to us.)
And thus the absence of a definition of order in a manner of speaking is responsible for the
subject of this paper: the Diophantine undecidability of arbitrary holomorphy rings of
characteristic 0 not equal to a field. The main results of this paper are stated below.
\begin{introthm}%
\label{thm:countable}
Let $K$ be a countable recursive one variable function field of characteristic 0.  Then Hilbert's Tenth Problem is not
solvable over any holomorphy ring of $K$ not equal to the whole field.
\end{introthm}%

\begin{introthm}%
\label{thm:uncountable}%
Let $K$ be a one variable function field of characteristic 0 over a field of constants $C$. Then for
any holomorphy ring of $K$ not equal to the whole field, there exist elements $x_1,\ldots, x_n \in
K \setminus C$ such that there is no algorithm to tell whether a polynomial equation with
coefficients in $\Q(x_1,\ldots,x_n)$ has solutions in the ring.
\end{introthm}%

\begin{introremark}%
The reason for two separate statements has to do with the possibility that the function field $K$ is
not countable.  That possibility forces us to examine more carefully what we can allow as
coefficients of our polynomial equations.  In the case of a countable field it is possible to
allow every element of the field as a coefficient, but in the case the field is not countable we
have to restrict the set of possible coefficients to a countable set.  In our case this set will
depend on the ring.
\end{introremark}%

\begin{introremark}%
\label{rem:S-integers}The case of holomorphy rings which are actually rings of $\calS$-integers, i.e. rings where only finitely
many primes of the field are allowed in the pole divisors of the ring elements, has been treated by the second author in
\cite{Sh92}. While it is not explicitly discussed in the paper, the statement of the Theorem \ref{thm:uncountable} follows from
the construction of the equations. We should also like to note that the aforementioned paper of the second author, just as the
present one, was a generalization of ideas of Denef from \cite{Den4}. In this paper Denef  used a Pell equation to construct a
model of integers instead of an elliptic curve. In view of this result we will always assume that the set of primes allowed in
the pole divisors is infinite.
\end{introremark}%

In some cases we will be able to prove a stronger result giving an existential definition of $\Z$
over an arbitrary holomorphy subring of the field not equal to the whole field. More specifically
the following theorem holds.

\begin{introthm}%
\label{thm:definable}
Let $K$ be any function field of characteristic 0 over a field of constants $C$.  Assume there
exists a subset $C_0$ of $C$ such that $C_0$ contains $\Z$ and has a Diophantine definition over
$K$.  Then $\Z$ is existentially definable over any holomorphy ring of $K$ not equal to the whole
field.
\end{introthm}%

\begin{introremark}%
We know of many function fields of characteristic 0 where constants are existentially definable.
They include function fields over ample fields of constants and other large fields, including
fields which are algebraically closed.  (See \cite{Koenig}, \cite{Pop}, and \cite{Sh35} for various
examples.)
\end{introremark}%

 The main idea behind the proofs of Theorems \ref{thm:countable} -- \ref{thm:definable} is rather
simple.  In a ring, where not all primes are inverted, there is a natural way to define the order
using divisibility. So even if we cannot do it over a field, we can define the order over a ring
(or come pretty close).  We are now ready to proceed with the technical details.
\section{Basic Diophantine Facts.}\label{section:diophantine}

We start with giving precise definitions to the objects we are going to study, beginning with Diophantine sets.

\begin{definition}%
\label{def:diophdefinition}
Let $R_0\subset R$ be rings and let $A \subset R^m$. A  \emph{Diophantine definition} of $A$ over $R$, with coefficients in $R_0$,
is a finite collection of polynomials $\{f_{i,j}(t_1,\ldots,t_m,x_1,\ldots,x_n), i=1,\ldots,r, j=1,\ldots,s\} \subset
R_0[t_1,\ldots,t_m,x_1,\ldots,x_n]$ such that for any $(t_1,\ldots,t_m) \in R^m$, we have the equivalence
\[ %
(t_1,\ldots,t_m)\in A\quad\Longleftrightarrow\quad
\exists x_1,\ldots,x_n \in R,
\bigvee_{j=1}^s\bigwedge_{i=1}^r f_{i,j}(t_1,\ldots,t_m,x_1,\dots,x_n) = 0.
\]%
We say that $A$ is \emph{Diophantine over $R$} w.r.t.~$R_0$ if it has such a Diophantine definition.
\end{definition}%
\begin{remark}
In general $R_0$ plays an auxiliary role and is often omitted, the default value being of course $R$.
\end{remark}
\begin{remark}
\label{rem:diophdefinition} Consider a Diophantine set $A$ as in  Definition \ref{def:diophdefinition}.

If $R$ is a domain (which is generally the case in applications), then $A$ has a Diophantine
definition ``with $s=1$'', i.e.~consisting of a system of polynomial equations, without the
disjunction operation.

If, moreover, the fraction field of $R$ is not algebraically closed, we can even take $r=s=1$; in
other words, $A$ has Diophantine definition consisting of one equation. Most authors take this as
the definition of a Diophantine set. (See \cite{Da2} or \cite{Sh34}, Chapter I for more details.)
\end{remark}

We will be able to construct such a definition of  $\Z$ over holomorphy subrings of our function
field $K$ provided a subset of the constant field containing $\Z$ has a Diophantine definition over
$K$.

\begin{definition}%
\label{def:H10}%
Let $R$ be a ring and let $R_0$ be a recursive subring of $R$. We say that \emph{Hilbert's tenth
problem is solvable} in $R$, with coefficients in $R_0$, if there is an algorithm taking as input a
finite set of polynomials in $R_0[X_1,\dots,X_m]$ (for some arbitrary $m>0$) and telling whether
they have a common zero in $R^m$.\\ We write H10($R, R_0$) for this property. If $R$ is
recursive we take H10($R$) to mean H10($R, R$).
\end{definition}%
\begin{remark}%
\label{rem:H10/1}%
Assume H10($R,R_0$) holds, and let $A\subset R^m$ be Diophantine w.r.t.~$R_0$, with given
Diophantine definition $(f_{i,j})_{1\leq i\leq r,1\leq j\leq s}$. Then $A$ is a finite union of
projections of sets $A_j$ ($1\leq j\leq s$) defined by polynomial systems (with some extra
variables). Since $A=\emptyset$ if and only if each $A_j$ is empty, there is an algorithm (taking
$(f_{i,j})$ as input) telling whether $A$ is empty or not. This of course could be taken as a
definition for the H10 property.

With the same assumptions, let $t$ be a point in $R_0^m\subset R^m$. Then $\{t\}$ is Diophantine
w.r.t.~$R_0$, and we have that $t\in A$ if and only if $\{t\}\cap A\neq\emptyset$. Hence, the above
discussion shows that there is an algorithm telling whether $t$ belongs to $A$. \end{remark}
\begin{remark}%
 \label{rem:H10/2} %
Just as in Remark \ref{rem:diophdefinition}, if $R$ is a domain with non-algebraically closed
fraction field, it suffices to check H10 for systems consisting of one polynomial: this is the
traditional definition of the H10 property.
\end{remark} %

\begin{proposition}
\label{prop:diophmodelcriterion}%
Let $R_1\subset R_2\subset R_3$ be rings, with $R_1$, $R_2$ and the inclusion $R_1\subset R_2$
recursive. Let $I$ be an ideal of $R_3$ with the following properties:
\begin{itemize}
\item $I$ is generated by finitely many elements of $R_2$ \emph{(in particular, it is Diophantine w.r.t.~$R_2$)}.
\item $R_1\cap I=\{0\}$.
\item The set $R_1+I\subset R_3$ is Diophantine w.r.t.~$R_2$.
\end{itemize}
Then \textup{H10($R_3,R_2$)} implies \textup{H10($R_1$)}.
\end{proposition}
\begin{proof} 
Assume \textup{H10($R_3,R_2$)}. Let $D\subset R_1^m$ be defined by
$$D:=\bigl\{ t\in R_1^m\mid \forall i\in\{1,\dots,r\},\;f_i(t)=0\bigr\}$$
where the $f_i$'s are polynomials with coefficients in $R_1$.  We are looking for an algorithm telling whether $D$ is empty.

Put $\Delta:=R_1+I$, which is Diophantine in $R_3$ by assumption, and define $B\subset R_3^m$ by
$$B:=\bigl\{t\in \Delta^m\,\mid \forall i\in\{1,\dots,r\},\;f_i(t)\in I\bigr\}.$$
Clearly, $B$ is Diophantine w.r.t.~$R_2$ since $\Delta$ and $I$ are. Hence, by \textup{H10($R_3,R_2$)}, there is an algorithm telling whether $B$ is empty, so it suffices to prove that $D=\emptyset$ if and only if $B=\emptyset$. The ``if'' part is trivial since $D\subset B$. Conversely, assume there exists some $t\in B$. By definition of $\Delta$, there exists $t_1\in R_1^m$ which is congruent (coordinatewise) to $t$ $\bmod\,I$. Then for each $i$ we still have $f_i(t_1)\in I$, but also $f_i(t_1)\in R_1$ since $f_i$ has coefficients in $R_1$. Hence $f_i(t_1)=0$, which means that $t_1\in D$, hence $D\neq\emptyset$. (In fact it is easy to see that $B=D+IR_3^m$ and $D=R_1^m\cap B$.)
\end{proof}

We will also use the following standard trick:

\begin{proposition}
\label{prop:freeextension}
Let $R_0\subset R\subset R'$ be rings, with $R_0$ recursive. Assume that, as an $R$-module, $R'$ has a finite basis $\mathcal{B}=\{b_1,\dots,b_m\}$ such that $R_0$ contains the following elements:
\begin{itemize}
\item the coordinates of $1$ in $\mathcal{B}$,
\item the entries of the matrix of multiplication by $b_i$ in $R'$, for each $i\in\{1,\dots,m\}$.
\end{itemize}
We identify $R'$ with $R^m$ using $\mathcal{B}$. 
Let $D\subset R'^d$ be Diophantine over $R'$ w.r.t.~$R_0$. Then:
\begin{enumerate}
\item\label{prop:freeextension1} $D$ (as a subset of  $R^{md}$) is Diophantine over $R$ (w.r.t.~$R_0$).
\item\label{prop:freeextension2} $D\cap R^d\subset R^d$ is Diophantine over $R$ (w.r.t.~$R_0$).
\end{enumerate}
\end{proposition}
\begin{proof} (\ref{prop:freeextension1}) is immediate from the assumptions. The first assumption also implies that the inclusion $R\subset R'\cong R^m$ identifies  $R$ with a Diophantine subset of $R^m$ w.r.t. $R_0$, so the same property holds for the inclusion $R^d\subset {R'}^d\cong R^{md}$. Assertion (\ref{prop:freeextension2}) follows.
\end{proof}

\section{Basic Facts on Function Fields and Holomorphy Rings.}\label{section:funfields}

\begin{definition}\label{def:funfield}%
Let $C$ be a field. Then a (one-variable) function field $K$ over $C$ is is a finite extension of the rational function field
$C(t)$. (Equivalently, it is a finitely generated extension of $C$ of transcendence degree $1$). For such a function field $K$, a
\emph{prime} of $K$ is a nontrivial discrete valuation of $K$ which is trivial on $C$. We denote by $\calP_K$ the set of such
primes (abusingly omitting $C$ from the notation). For $\pp\in\calP_K$ we adopt the traditional notation $\ord_\pp$ for the
corresponding normalized valuation and we denote by $O(\pp)$ the associated valuation ring. If $\calW$ is a \emph{non-empty}
subset of $\calP_K$, we put

\[%
O_{K,\calW}=\bigcap_{\pp\not\in\calW}O(\pp)= \bigl\{h \in K\mid
\forall \pp \not \in \calW \mbox{ we have that } \ord_{\pp}h
\geq 0\bigr\}.
\]%
$O_{K,\calW}$ is called a \emph{holomorphy ring} of $K$.
\end{definition}%

Note that, with the above notations,  taking $\calW=\emptyset$ would lead to the intersection of
\emph{all} rings $O(\pp)$ ($\pp\in\calP_K$). This ring is the algebraic
closure $C'$ of $C$ in $K$, a finite extension of $C$,
and $K$ is a function field over $C'$ with the same set
$\calP_K$ and rings $O_{K,\calW}$ as over $C$. Therefore, for
the purposes of this paper we can always replace $C$ by $C'$,
i.e.~ assume $C$ algebraically closed in $K$. (In other words,
$K$ is a \emph{regular extension} of $C$).

If $\calW=\calP_K$, then $O_{K,\calW}=K$. Otherwise, $O_{K,\calW}$ is a
Dedekind domain with fraction field $K$; its maximal ideals correspond bijectively to elements of $\calP_K\setminus\calW$, via
the map
\[%
\pp\longmapsto I_\pp:=\bigl\{h\in O_{K,\calW}\mid \ord_{\pp}h > 0\bigr\}.
\]%
We shall derive Theorems \ref{thm:countable}, \ref{thm:uncountable} and \ref{thm:definable} from the following result:
\begin{theorem}\label{theorem:diophantine} Let $K$ and $\calW$ be as above, with $\emptyset\neq\calW\neq\calP_K$ and $\mathrm{char}\,K=0$.  
Let $\pp\not\in\calW$ be a prime of $K$, and let $I_\pp$ be the corresponding maximal ideal of $O_{K,\calW}$. Then $\Z+I_\pp$ is Diophantine in 
 $O_{K,\calW}$.
 \end{theorem}
\noindent\emph{Proof of Theorems \ref{thm:countable}, \ref{thm:uncountable} and \ref{thm:definable}\/} (from \ref{theorem:diophantine}): For Theorems \ref{thm:countable} and \ref{thm:uncountable} we simply apply Proposition
\ref{prop:diophmodelcriterion} with $R_1=\Z$ and $R_3=O_{K,\calW}$. We  take for $R_2$ the subring
of $R_3$ generated by a finite subset containing a generating set for $I_\pp$ and all elements occurring in a Diophantine definition of $\Z+I_\pp$. The condition $I_\pp\cap\Z=\{0\}$ is obvious since nonzero integers are invertible in
$C\subset O_{K,\calW}$.

For Theorem \ref{thm:definable}, we have $\Z\subset(\Z+ I_\pp)\cap C_0\subset(\Z+ I_\pp)\cap C=\Z$ (the last equality follows from $C\cap I_\pp=\{0\}$). Hence $\Z=(\Z+ I_\pp)\cap C_0$ is Diophantine.\qed
\medskip
 
From now on, we assume that $C$ is algebraically closed in $K$, and that the characteristic is zero. If $\hat{C}$ is a finite
extension of $C$, then $\hat{K}:=\hat{C}\otimes_C K$ is a finite extension of $K$, and a function field over $\hat{C}$. Moreover,
for $\calW$ as above, the following three subrings of $\hat{K}$ are equal:%
\begin{itemize}%
\item $\hat{C}\otimes_C O_{K,\calW}$,%
\item the integral closure of $O_{K,\calW}$ in $\hat{K}$,%
\item the holomorphy ring $O_{\hat{K},\hat{\calW}}$ where $\hat{\calW}$ is the set of primes of $\hat{K}$ inducing primes in
$\calW$.%
\end{itemize}%
The first description shows in particular that $O_{\hat{K},\hat{\calW}}$ is a free module over $O_{K,\calW}$, of rank
$[\hat{C}:C]$. 

In addition, let $\pp\not\in\calW$ be a prime of $K$. 
Then there exists a prime $\hat{\pp}\not\in\hat{\calW}$ of $\hat{K}$ extending $\pp$, and the corresponding 
ideals $I_\pp\subset O_{K,\calW}$ and $I_{\hat{\pp}}\subset O_{\hat{K},\hat{\calW}}$ satisfy $I_\pp=O_{K,\calW}\cap I_{\hat{\pp}}$. 
It follows that $\Z+I_\pp=O_{K,\calW}\cap (\Z+I_{\hat{\pp}})$. By Proposition \ref{prop:freeextension}\,(\ref{prop:freeextension2}) we see 
that if $\Z+I_{\hat{\pp}}$ is Diophantine in $O_{\hat{K},\hat{\calW}}$, then $\Z+I_{{\pp}}$ is Diophantine in $O_{{K},{\calW}}$. 

In particular, to prove Theorem \ref{theorem:diophantine} for $K$, $\calW$, $\pp$, we may replace these data by  $\hat{K}$, $\hat{\calW}$, $\hat{\pp}$, respectively.

We shall use this remark as follows: take for $\hat{C}$ the residue field of $\pp$. We have a
surjective morphism $O_{K,\calW}\to\hat{C}$ of $C$-algebras, hence (tensoring with $\hat{C}$) a surjective $\hat{C}$-morphism
$O_{\hat{K},\hat{\calW}}\to\hat{C}$. Its kernel defines a prime $\hat{\pp}$ of $\hat{K}$ above $\pp$, which has \emph{degree one} in the sense that its residue field is the constant field $\hat{C}$ of $\hat{K}$. To summarize,
\emph{we can always assume that the prime $\pp$ of Theorem  \ref{theorem:diophantine} has degree one. }

As we have already mentioned above,  our paper has two main inputs.  The first one is contained in a
paper of  Denef (see \cite{Den4}) which  constructs a rank one elliptic curve over any rational
function field of characteristic 0 together with a way of generating integers.  The second input is
a  result of the first author that allows the elliptic curve constructed by Denef to retain its nice
properties under finite extensions.  More specifically we will use the following result which is a
 consequence of Theorem 1.8 (ii) and Proposition 2.3.1 of \cite{MB3}.
\begin{theorem}%
\label{thm:MB} %
Let $K$ be a function field of characteristic 0 over a field of constants $C$. Let $\pp$ be a
degree one prime of $K$. Let $\Dd$ be a divisor of $K$ such that $\ord_{\qq}\Dd \in\{0,1\}$ 
for any prime $\qq$ of $K$, $\ord_{\pp}\Dd=0$,
and the degree of $\Dd$ is at least $2g_K+2$, where $g_K$ is the genus of $K$. Let
$F(T)$ be a nonsingular cubic polynomial over $\Q$ such that the elliptic curve $Y^2=F(X)$ has no
complex multiplication. Then there exists an $x \in K$ such that its pole divisor is $\Dd$,
$\ord_{\pp}x >0$, and the elliptic curve $E_x$ defined by the equation
\begin{equation}%
\label{eq:MD}%
F\Bigl(\frac{1}{x}\Bigr)\,Y^2=F(X)
\end{equation}%
has the property that $E_x(C(x)) = E_x(K)$. Also $E_x(C(x))$ is of rank 1 generated
by the point with affine coordinates $(\frac{1}{x},1)  \in E(C(x)) \setminus E(C)$ modulo 2-torsion.
\end{theorem}%
\begin{proof}%
We need a slight refinement of   Proposition 2.3.1 of \cite{MB3}:
\begin{lemma}%
\label{le:function}
With the assumptions of Theorem \ref{thm:MB}, there exists a nonzero $g\in K$ with the following properties:
\begin{itemize}
\item the divisor of zeros of $g$ is $\Dd$,
\item $g$ has only simple poles, and $\pp$ is one of them,
\item $g$ has  simple ramification (i.e.~in the extension $K/C(g)$ no prime has ramification degree greater than 2).
\end{itemize}
\end{lemma}
\begin{proof}%
The argument is classical and entirely similar to \cite{MB3}, 2.3.1. Put $d:=\deg\Dd$. The
linear system $\vert\Dd\vert$ is a projective space of dimension $d-g$, and we identify
$\Dd$ with a ($C$-rational) point in it. Inside $\vert\Dd\vert$ we consider the following
subvarieties, where $Q$ (resp.~$\Delta$) denotes a variable point (resp.~effective divisor):
\begin{itemize}%
\item $H=\{\text{divisors of the form } \pp+\Delta\}$,%
\item $Z_1=\{\text{divisors of the form } 2\pp+\Delta\}$,%
\item $Z_2=\{\text{divisors of the form } \pp+2Q+\Delta\}$.%
 \item $Z_3=\{\text{divisors of the form } 3Q+\Delta\}$,%
 \end{itemize}%
Clearly, none of these contains $\Dd$, and $H$ is a hyperplane because $\pp$ has degree $1$. It
is proved in \cite{MB3} that $Z_3$ has codimension $\geq2$ in $\vert\Dd\vert$, and similar
arguments easily show that the same holds for $Z_1$ and $Z_2$. Hence we can find a line in
$\vert\Dd\vert$ through the point $\Dd$, defined over $C$ and disjoint from $Z_1\cup
Z_2\cup Z_3$. This line meets $H$ at a point $\Dd'$. There is an element $g$ of $K$ with
divisor $\Dd-\Dd'$, and this $g$ has the required properties.
\end{proof}%

Let us return to the proof of Theorem \ref{thm:MB}. Clearly, by multiplying $g$ by some nonzero
constant $c$ (in $\Q$, if we wish) we may also choose $g$ with branch locus disjoint from the
inverse roots of $F$. This makes $g$ \emph{admissible} for $\Dd$ in the sense of \cite{MB3},
Definition 1.5.2. Further, it follows from Theorem 1.8 (ii) of \cite{MB3} that by choosing $c$
appropriately we may assume in addition that $g$ is \emph{good}, i.e.~$E_{g^{-1}}(C(g))
=E_{g^{-1}}(K)$.

Now let $x = g^{-1}$: we now have that $E_x(C(x)) = E_x(K)$, and the rest of the theorem follows
from the assumption on $F$ and from \cite{Den4} which describes $E_x(C(x))$.
\end{proof}%

\section{Diophantine Undecidability of Holomorphy Rings.}\label{section:UndecHolrings}

\begin{notationassumption}%
\label{not:1}
We start with a first notation set.
\begin{itemize}
\item Let $K$ be a function field of characteristic 0 over a field of constants $C$. %
\item Let $\calP_K$ be the set of all primes of $K$.%
\item Let $\emptyset \not =\calW \subset \calP$, $\calW \not = \calP$.%
\item Let $\pp \in \calP \setminus \calW$ be a prime of degree $1$. (By assumption, $\calP \setminus \calW$ is not empty, and as
explained in the previous section we may assume that it contains a prime of degree $1$ by extending $K$).%

\item Let $g_K$ be the genus of $K$. %
\item Assume that $\calW$ contains infinitely many primes (in fact the present proof works whenever $\sum_{\qq\in\calW}\deg\qq\geq 2g_K+2$). As noted in Remark \ref{rem:S-integers} of the introduction, the case where $\calW$ is finite is settled in \cite{Sh92}; more precisely, in that case, $\Z$ is Diophantine in  $O_{K,\calW}$ (Theorem 3.1 of \cite{Sh92}), hence our Theorem \ref{theorem:diophantine} also holds.
\item Let $\calD=\{\qq_1, \ldots, \qq_{2g_K+2}\}$ be a set of distinct elements of $\calW$. (We
only need the total degree of $\calD$ to be at least $2g_K+2$, in order to apply Theorem
\ref{thm:MB}.)

\item Let $x\in K$ be such that its pole divisor is $\prod_i \qq_i$, $\ord_{\pp}x>0$, and $E(K)=E(C(x))$, where $E=E_x$ is the
elliptic curve defined in (\ref{eq:MD}). (Such an $x$ exists by Theorem \ref{thm:MB}).
 \end{itemize}%
\end{notationassumption}%

Observe that $x \in O_{K,\calW}$ and therefore $O_{K,\calW}$ contains the polynomial ring $C[x]$. Moreover, the condition
$\ord_{\pp}x>0$ means that $\pp$ lies above the ideal $(x)$ of $C[x]$. In other words, for any $z\in C(x)$ we have
\begin{equation}\label{equation:order}
\ord_{\pp}z=(\ord_{\pp}x)(\ord_{0}z)
\end{equation}
where, in the right-hand side, $z$  is viewed as a rational function of $x$.\\

Two elements $a$, $b \in O_{K,\calW}$ will be called \emph{coprime} if they generate the unit ideal, i.e.\ there
exist $A, B \in O_{K, \calW}$ such that
\begin{equation}%
\label{eq:relprime}%
Aa +Bb = 1.
\end{equation}%
Note that we have ``Gauss' Lemma'': if $a$ and $b$ are coprime and $a$ divides $bc$ in $O_{K,\calW}$, then $a$ divides $c$.

\begin{proposition}%
The set $\bigl\{h \in O_{K,\calW}\mid h \not=0\bigr\}$ is Diophantine over $O_{K,\calW}$.
\end{proposition}%
\begin{proof}%
The proof is essentially a consequence of the Strong Approximation Theorem and can be found in \cite{Sh5}.
\end{proof}%

\begin{notationassumption}%
\label{not:2}
We now add the following notation and assumptions to our list.
\begin{itemize}%
\item Let $P \in E(K)$ be the point whose affine coordinates derived from (\ref{eq:MD}) are $(\frac{1}{x}, 1)$.%
\item For nonzero $n\in\Z$, let $(x_n,y_n)$ be the affine coordinates of $[n]P$ derived from (\ref{eq:MD}). Since $P\in E(C(x))$
we have that $x_n$ and $y_n$ are rational functions of $x$. \item Since $[n]P$ is not
a torsion point, we have $y_n\neq0$ and we can write
\[%
\frac{xx_n}{y_n}=\frac{\alpha_n}{\beta_n},
\]%
where $\alpha_n, \beta_n \in C[x]$ are relatively prime polynomials in $x$ (in
particular, since they satisfy relation (\ref{eq:relprime}) in $C[x]\subset O_{K,\calW}$, they are also
coprime in $O_{K,\calW}$).
\end{itemize}%
\end{notationassumption}%
The following lemma shows how we will generate integers to show undecidability. Its proof can be
found in Lemma 3.2 of \cite{Den4}.
\begin{lemma}%
\label{le:equiv}%
For any $n \in \Z_{>0}$ we have that $\ord_{\pp}\bigl(x\frac{ x_n}{y_n}-n\bigr)>0$.\qed%
\end{lemma}%
Using Lemma \ref{le:equiv} and the definition of being coprime, one easily deduces the following lemma.
\begin{lemma}%
\label{le:ord0}
Let $n$ be a nonzero integer. Assume that
$$\frac{xx_n}{y_n}=\frac{a_n}{b_n}$$
where $a_n$, $b_n\in O_{K,\calW}$ are coprime.  Then:
\begin{eqnarray}%
\label{le:ord0;eq:0}
a_n=\varepsilon\alpha_n\quad\textup{ and } \quad b_n=\makebox[2em][l]{$\varepsilon\beta_n$}&  &\textup{for some }\varepsilon \in O_{K,\calW}^\times\,,\\
\label{le:ord0;eq:1}
a_n-n\,b_n=\makebox[2em][l]{$xw$} && \textup{for some }w \in O_{K,\calW}\,,\\
\label{le:ord0;eq:2}
\ord_\pp(b_n)=\makebox[2em][l]{$0.$}&&
\end{eqnarray}
\end{lemma}%
\begin{proof}%
To prove (\ref{le:ord0;eq:0}) we note that $a_n\beta_n=b_n\alpha_n$, hence (by Gauss' lemma) $a_n$ and $\alpha_n$ divide each other in $O_{K,\calW}$. In other words, $\varepsilon:=a_n/\alpha_n$ is a unit.

To prove (\ref{le:ord0;eq:1}) and (\ref{le:ord0;eq:2}) we may and will assume, in view of (\ref{le:ord0;eq:0}), that
$a_n=\alpha_n$ and $b_n=\beta_n$. Let $w\in K$ be defined by (\ref{le:ord0;eq:1}). Let us prove that $w\in C[x]$. Since $a_n$
and $b_n$ are in $C[x]$ it suffices to prove that $a_n-n\,b_n$ (as a polynomial in $x$) vanishes at $0$, which by
(\ref{equation:order}) is equivalent to $\ord_\pp (a_n-n\,b_n)>0$. This is clear from Lemma \ref{le:equiv} since
$a_n-n\,b_n=b_n(x\frac{ x_n}{y_n}-n)$ and $b_n\in C[x]$.

Let us now prove (\ref{le:ord0;eq:2}). Since $a_n-n\,b_n$ vanishes at $0$ it follows that if $b_n$
vanishes at $0$, so does $a_n$. Since they are relatively prime polynomials, this cannot happen.
\end{proof}%

We also have a converse of sorts to (\ref{le:ord0;eq:1}) above.
\begin{lemma}
\label{le:ord1}%
With the assumptions of \ref{le:ord0}, suppose for some $c \in K$ we have that $\ord_{\pp}(a_n-b_nc)>0$. Then $\ord_{\pp}(c
-n)>0$.
\end{lemma}%
\begin{proof}%
The equation $a_n-b_nc=xw$ implies that $\ord_{\pp}(a_n-b_nc) >0$ as $\ord_{\pp}x >0$ and $\pp \not
\in \calW$. Since $\ord_{\pp}(a_n-b_nn)>0$ by Lemma \ref{le:ord0}, we conclude that
$\ord_{\pp}(b_n(c-n))>0$. This proves the result by (\ref{le:ord0;eq:2}).
\end{proof}%
Next we prove an easy lemma.
\begin{lemma}
\label{le:Edioph}%
The set
\[%
{\tt E}=\left\{(u,v,w,z)\in (O_{K,\calW})^4\;\Big{|}\; vw\neq0,\;\exists n \in \Z \setminus \{0\}: x_n=\frac{u}{v},
y_n= \frac{z}{w}\right\}
\]%
is Diophantine.
\end{lemma}%
\begin{proof}%
Since we know how to define non-zero elements over any holomorphy ring and all the points of $E(K)$
are in fact of the form $[n]P + T$ where $T$ is a 2-torsion point, we can easily define the set
\[%
{\tt E}_{\mbox{even}}=\left\{(u',v',w',z')\in (O_{K,\calW})^4\,\Big{|} \,\exists k \in \Z \setminus \{0\}:
x_{2k}=\frac{u'}{v'}, y_{2k}= \frac{z'}{w'}\right\}
\]%
using Theorem \ref{thm:MB}.  Then $(u,v,w,z) \in {\tt E}$ if and only if either $(u,v,w,z) \in {\tt
E}_{\mbox{even}}$ or $\bigl(\frac{u}{v}, \frac{z}{w}\bigr) = \bigl(\frac{u'}{v'},\frac{z'}{w'}\bigr) +_E
\bigl(\frac{1}{x},1 \bigr)$, where $ (u',v',w',z') \in {\tt E}_{\mbox{even}}$.
\end{proof}%

\noindent\emph{Proof of Theorem  \ref{theorem:diophantine}:\/} Put $I=I_\pp=\left\{t\in O_{K,\calW}\mid \ord_\pp(t)>0\right\}$, and let us prove that $\Z+I$ is Diophantine in $O_{K,\calW}$.

 Let $\xi$ be an element of $O_{K,\calW}$. We claim that the following are equivalent:
\begin{enumerate}
  \item\label{prop:Zdioph1} $\xi\in\Z+I$,
  \item\label{prop:Zdioph2} either $\xi\in I$, or the following system has a solution $(u,v,w,z,a,b,A,B)$ in $O_{K,\calW}^8$:
  \begin{equation}%
\label{eq:sys2}
\left \{
\begin{array}{c}
\rule[-3mm]{0mm}{0mm}(u,v, w,z) \in {\tt E}\\
\displaystyle \frac{\,a\,}{b}=\frac{\,xuw\,}{vz}\\
\rule[-3mm]{0mm}{8mm}Aa + Bb=1\\
a-b\xi\in I.
\end{array}
\right .
\end{equation}%
\end{enumerate}

This clearly implies the result since both ${\tt E}$ and $I$ are Diophantine (the
former by \ref{le:Edioph}, and the latter because it is a finitely generated ideal).

First, assume (\ref{prop:Zdioph1}). If $\xi\in I$ we are done. Otherwise, we may assume that $\xi=:n$ is a nonzero integer since
both (\ref{prop:Zdioph1}) and (\ref{prop:Zdioph2}) are invariant under adding an element of $ I$ to $\xi$. We construct a
solution of (\ref{eq:sys2}) as follows. First, choose $u, v, w,z$ so that $\frac{u}{v}=x_{n}, \frac{z}{w} = y_{n}$
(first relation). Put $a=\alpha_{n}$ and $b=\beta_{n}$ (as defined in \ref{not:2}): the second relation is satisified and we can find $A$ and $B$ satisfying
the third. Finally by (\ref{le:ord0;eq:1}) the fourth relation holds since $x\in I$.

Now assume that (\ref{prop:Zdioph2}) holds. As before, (\ref{prop:Zdioph1}) is trivial if $\xi\in I$. Otherwise, fix a solution
$(u,v,w,z,a,b,A,B)$ of (\ref{eq:sys2}). By definition of $\tt E$ there is a nonzero integer $n$ such that $x_n=\frac{u}{v}$ and
$y_n= \frac{z}{w}$. We can then apply Lemmas \ref{le:ord0} and \ref{le:ord1} with $a_n=a$, $b_n=b$ and $c=\xi$ to conclude that
$\ord_{\pp}(\xi-n) >0$. In other words, $\xi\in\Z+ I$. \qed

\end{document}